\documentclass[11pt,a4paper]{article}
\usepackage[utf8x]{inputenc}
\usepackage[english]{babel}
\usepackage[T1]{fontenc}
\usepackage{amsmath}
\usepackage{amsfonts}
\usepackage{amssymb}
\usepackage{makeidx}
\usepackage{graphicx}
\usepackage[left=3cm,right=2cm,top=3cm,bottom=2cm]{geometry}
\usepackage{setspace}
\usepackage{amsthm}
\usepackage{cite}

\usepackage[colorlinks=true,urlcolor=blue,
citecolor=red,linkcolor=blue,linktocpage,pdfpagelabels,
bookmarksnumbered,bookmarksopen]{hyperref}

\addto\captionsenglish{}

\newtheorem{definition}{Definition}[section]

\newtheorem{proposition}{Proposition}[section]

\newtheorem{theorem}{Theorem}[section]

\newtheorem{example}{Example}[section]

\newtheorem{lemma}{Lemma}[section]

\newtheorem{observation}{Remark}[section]

\newtheorem{corollary}{Corolary}[section]

\numberwithin{equation}{section}

\newcommand{\bo}{\begin{observation}}
\newcommand{\eo}{\end{observation}}
\newcommand{\bd}{\begin{definition}}
\newcommand{\ed}{\end{definition}}
\newcommand{\bp}{\begin{proposition}}
\newcommand{\ep}{\end{proposition}}
\newcommand{\bt}{\begin{theorem}}
\newcommand{\et}{\end{theorem}}
\newcommand{\bc}{\begin{corollary}}
\newcommand{\ec}{\end{corollary}}
\newcommand{\bl}{\begin{lemma}}
\newcommand{\el}{\end{lemma}}
\newcommand{\be}{\begin{example}}
\newcommand{\ee}{\end{example}}
\newcommand{\beq}{\begin{equation}}
\newcommand{\eeq}{\end{equation}}
\newcommand{\beqa}{\begin{equation*}}
\newcommand{\eeqa}{\end{equation*}}
\newcommand{\R}{\mathbb{R}}
\newcommand{\RN}{\mathbb{R}^{N}}
\newcommand{\N}{\mathbb{N}}
\newcommand{\WumN}{ W^{1,N}(\mathbb{R}^{N}) }

\newcommand{\LN}{ L^{N}(\mathbb{R}^N) }

\newcommand{\LNquatro}{L^{\frac{2N^2}{2N-1}}(\RN)}

\newcommand{\Nn}{\mathcal{N}}

\newcommand{\Ls}{ L^{s}(\mathbb{R}^N)}

\newcommand{\intR}{\displaystyle\int\limits_{\mathbb{R}^N}}
\newcommand{\un}{u_{n}}

\newcommand{\uk}{u_{k}}

\newcommand{\vn}{v_{n}}

\newcommand{\yn}{y_{n}}

\newcommand{\until}{\tilde{u}_n}

\newcommand{\RA}{\rightarrow}

\usepackage{times, color, xcolor}
\addto\captionsportuguese{
}

\begin{document}

 	\title{On the existence and multiplicity of solutions for the $ N $-Choquard logarithmic equation with exponential critical growth
		}
	\author{
		Eduardo  de S. Böer \thanks{Corresponding author} \thanks{ E-mail address: eduardoboer04@gmail.com Tel. +55.51.993673377, https://orcid.org/0000-0002-3401-3702.}  and Ol\'{\i}mpio H. Miyagaki  \footnote{ E-mail address: ohmiyagaki@gmail.com, Tel.: +55.16.33519178 (UFSCar), https://orcid.org/0000-0002-5608-3760.}\\
		{\footnotesize Department of Mathematics, Federal University of S\~ao Carlos,}\\
		{\footnotesize 13565-905 S\~ao Carlos, SP - Brazil}\\ }
\noindent
				
	\maketitle

\noindent \textbf{Abstract:} The present work is concerned with the Choquard Logarithmic equation \linebreak $-\Delta_N u + a(x)|u|^{N-2}u + \lambda (\ln|\cdot|\ast |u|^{N})|u|^{N-2}u = f(u) \textrm{ \ in \ } \mathbb{R}^N,$ where $ a:\RN \RA \R $ is a continuous non-negative function, $ \lambda >0 $, $ N \geq 3 $ and $f: \R \RA [0, \infty) $ is continuous function that behaves like $ \exp(\alpha |u|^{\frac{N}{N-1}}) $ at infinity, for $ \alpha >0 $. We prove the existence of a nontrivial solution at the mountain pass level and a nontrivial ground state solution. Also, a multiplicity result is provided for the subcritical case, using genus theory.

\vspace{0.5 cm}


\section{Introduction}

In the present paper we are interested with the existence of solutions for the following class of equations
\begin{equation}\label{i1}
-\Delta_p u + a(x)|u|^{p-2}u + \lambda \Psi |u|^{p-2}u = 0 \textrm{ \ in \ } \mathbb{R}^N,
\end{equation}
where $ a:\RN \RA \R $ is a continuous non-negative function, $ \lambda >0 $, $ N \geq 3 $ and $ \Psi : \RN \RA \R $ stands for an internal potential for a nonlocal self-interaction of $ u $. 

Consider $ u $ a radial function, that is, $u(x)=u(r),\ r=|x|,$ $ r>0 .$  Then, in similarity with the case $ p=2 $, we have a \textit{``fundamental solution for the $ p $-laplacian''} equation, $ -\Delta_p u =0 $, given by
$$
\Gamma_p (x) = \left\{ \begin{array}{rcl}
\dfrac{1}{(p-N)\omega_{N-1}}|x|^{1-\frac{N-1}{p-1}} &\textrm{ if \ } &N\neq p ,\\
\dfrac{1}{\omega_{N-1}}\ln |x|& \textrm{ if \ }& N= p ,
\end{array} \right.
$$ 
where $ \omega_{N-1} $ is the $ (N-1)$-dimensional measure of the $ (N-1)$-sphere.

In this sense, we can consider equation (\ref{i1}) with $ \Psi = \Gamma_p \ast |u|^p $, that is
\begin{equation}
-\Delta_p u + a(x)|u|^{p-2}u + \lambda (\Gamma_p \ast |u|^p) |u|^{p-2}u = 0 \textrm{ \ in \ } \mathbb{R}^N .
\end{equation}
In this work we will focus on finding solutions for the p-planar case, $ N=p $, and with a nonlinearity of Moser-Trudinger type,  
\begin{equation} \label{P}
-\Delta_N u + a(x)|u|^{N-2}u + \lambda (\ln|\cdot|\ast |u|^{N})|u|^{N-2}u = f(u) \textrm{ \ in \ } \mathbb{R}^N,
\end{equation}
with $ a \equiv 1 $, $ \lambda = \omega_{N-1} $, $ N \geq 3 $ and $f: \R \RA [0, \infty) $ a continuous function that behaves like $ \exp(\alpha |u|^{\frac{N}{N-1}}) $ at infinity, with primitive $ F(s)=\int\limits_{0}^{s}f(t)dt $.

In the sequence, we make a quick overview of literature. We start mentioning that Choquard equations are well-known for its great range of applications in science, specially in physics and mechanics.  The case $ N=3 $ for Laplacian operator has been extensively studied, due to its relevance in physics. Although this equation is called ``Choquard equation'', its seminal paper is due to Fröhlich and Pekar in \cite{[12] , [11] , [22]}, where the authors describe the quantum mechanics of a polaron at rest, in the particular case, when $ V(x)\equiv a > 0 $ and $ \gamma > 0 $. Then, it was considered by Choquard in 1976, while studying an electron trapped in its hole. We highlight that the local nonlinear terms on the right side of these equations, such as $ f(u)= b|u|^{p-2}u $, for $ b\in \R $ and $ p>2 $, stands, usually, in Schrödinger equations as a modelling of the interaction among particles. 

Since the mentioned seminal works, several variations of Choquard equations have been considered as a way to model different phenomena. As an additional example, we cite \cite{[18]} where Penrose has derived a Choquard equation while discussing about the self gravitational collapse of a quantum-mechanical system. Moreover, many results have been derived recently about the existence and regularity of solutions, including more general convolution potentials, such as in the evolution equation $ i\partial_t \phi = \Delta \phi +(V \ast |\phi|^2)\phi,  $ which models the interaction of a large system of nonrelativistic bosonic atoms and molecules.  

We also refer the reader to \cite{[20] , [16]}, in which the authors treat variations of Choquard equation with $ N\geq 3 $ and $ p=2 $. 

When dealing with specifically logarithmic forms of Choquard equations, one has scantier literature. It took some time for scientists to figure out how to treat this case, because, the presence of the logarithmic function rises a bunch of difficulties. Specially while dealing with Cerami sequences. In this sense, we can cite some seminal and recent works \cite{[6], [cjj], [10], [4]} where the authors have proved the existence and multiplicity of solutions for Laplacian equations with polynomial nonlinearity. Moreover, in \cite{[alves]}, the authors have dealt with equation (\ref{P}) considering a nonlinearity with exponential critical growth and Laplacian operator. Finally, in \cite{[boer2]}, the authors prove existence and multiplicity results for the fractional $ p-$Laplacian operator with exponential growth nonlinearity. In problems with such critical behaviour, some additional difficulties arises while dealing with Cerami sequences and that is why authors need some stronger conditions over the nonlinearity. 

On the other hand, problems involving the $N$-Laplacian operator are important in many fields of sciences. Since they can accurately describe the behaviour of electric, gravitational and fluid potentials, they mostly appear in the fields of electromagnetism, astronomy and fluid dynamics. Moreover, the study of ground state solutions for the $N$-Laplacian equation is crucial in the research of evolutions equations which appear in non-Newton fluids, turbulent flows in porous media and other contexts. 

We refer the reader to \cite{[medeiros], [lam], [zhang]} and the references therein, in which the authors work or refer to work where many  different  cases of equations involving the $N$-Laplacian operator are studied, such as  bounded  domains  and  unbounded  domains,  different behaviour of nonlinearities, different types of boundary conditions, and so on. Particularly, many works focussed on the study of subcritical and critical growth for the involved nonlinearities which allows us to treat the problem variationally.

Our work intend to extend or complement the results found in the above papers. The main features of the class of problems considered in this paper, are that it is defined in the whole $\RN$, it has a convolution term with the logarithmic function, which changes sign and in unbounded from above and bellow, involves an exponential critical growth and the nonlinear operator $N$-Laplacian. Such combination of element drives to new possibilities of applications. 

Now we present the conditions imposed over the nonlinearity $ f $ in order to obtain the desired results. Such hypothesis are usual in works involving a Moser-Trudinger inequality and the $ N $-Laplacian operator, as one can see in \cite{[doO], [zhang], [lam]}.

In order to simplify some calculations and make the notation more concise, we will write
$$ R(\alpha , s) = \exp(\alpha |s|^{\frac{N}{N-1}}) - S_{N-2}(\alpha , s) = \sum\limits_{N-1}^{+\infty} \dfrac{\alpha^k}{k!} |s|^{\frac{N}{N-1}k} , $$
where $S_{N-2}(\alpha_0 , s) = \sum\limits_{k=0}^{N-2}\dfrac{\alpha_{0}^{k}}{k!}|s|^{\frac{N}{N-1}k}$. 

We recall that a function $ h $ is said to have \textit{subcritical} exponential growth at $ +\infty $, if
$$
\lim\limits_{s\RA + \infty}\dfrac{h(s)}{R(\alpha , s)} = 0 \textrm{ \ , for all \ } \alpha >0 ,
$$
and we say that $ h $ has $ \alpha_0 $-\textit{critical} exponential growth at $ +\infty $, if
$$
\lim\limits_{s\RA + \infty}\dfrac{h(s)}{R(\alpha , s)} = \left\{ \begin{array}{ll}
0, \ \ \ \forall \ \alpha > \alpha_{0} \\
+\infty , \ \ \ \forall \ \alpha < \alpha_0
\end{array} \right. .
$$
In this sense, we assume that $f$ has critical growth. Moreover, motivated by Lemma \ref{l1}, we assume the following growth condition
$$ f\in C(\R , \R), f(0)=0 \mbox{ and has critical exponential growth.} \leqno{(f_1)} $$
$$ \lim\limits_{s\RA 0^{+}} \dfrac{f(s)}{|s|^{N-1}}=0.  \leqno{(f_2)}$$  

Before giving the next condition, we point out that, following the idea introduced by Stubbe \cite{[21]}, we will consider the slightly smaller Hilbert space
\beq \label{X}
X = \left\{ u\in \WumN  ; \intR \ln(1+|x|)|u(x)|^N dx < \infty \right\} \subset \WumN ,
\eeq
endowed with the norm
$$
||u||_{X}^{N}=||u||^{N}+||u||_{\ast}^{N} \mbox{ \ \ , \ \ where \ \ } ||u||_{\ast}^{N}=\intR \ln(1+|x|)|u(x)|^{N} dx , 
$$
in which the portion of the associated functional involving the $ \ln |\cdot| $ function will be finite.

Next, we need some kind of Ambrosetti-Rabinowitz condition to guarantee that Cerami sequences in $ X $ are bounded in $ \WumN $. The reader can see that this condition is a bit stronger that some usually found in the literature. It happens because two main causes, which are the integral term involving the logarithmic and the exponential behaviour of $ f $ at infinity. Once we need a bit more information about $ ||\nabla \cdot||_N $ to get the boundedness of the integral involving the exponential term, we cannot even adapt the conditions of \cite{[6], [cjj], [10]}, for example. In this sense, we ask for

$$ \mbox{ there exists }\ \theta \geq  2N \ \mbox{ such that}\   f(s)s \geq \theta F(s) > 0, \ \mbox{for all }\  s\in \R \setminus \{0\}. \leqno{(f_3)}$$

However, while the above condition exclude the necessity of dealing with the integral involving the exponential term, it is not always possible to argue that way. So, based on works such as \cite{[9]}, we insert a condition that makes possible to get a boundedness for the integral involving the exponential term and add some other geometric properties for the associated functional. 

$$ \mbox{ there is}\  q>2N \ \mbox{ and}\  C_q> \dfrac{S_{q}^{q}[2(q-N)]^{}\frac{q-N}{N}}{\rho_{0}^{q-N}q^{\frac{q}{N}}} \ \mbox{ such that}\  F(s) \geq C_q |s|^q , \ \mbox{for all}\  s\in \R , \leqno{(f_4)}$$ 
where $ S_q $ is a constant obtained from the Sobolev embeddings and $ \rho_0 $ is a suitable value such that the exponential term can be controlled, defined in Lemma \ref{l10}.

Our study is based on variational techniques. In this sense, we will give some necessary informations about the problem (\ref{P}). One can see that the associated functional $ I: \WumN \RA \R \cup \{ \infty\} $ is given by
$$
I(u)=\dfrac{1}{N}\intR |\nabla u|^N + |u|^N dx + \dfrac{1}{2N} \intR \intR \ln(1+|x-y|) |u(x)|^N |u(y)|^N dx dy - \intR F(u) dx .
$$ 
As we mentioned earlier, $ I $ is not well-defined in $ \WumN $. But, it is possible to see that, for any $ u\in X $, from Moser-Trudinger Lemma \ref{l1}, $ I(u)< +\infty $. Moreover, by standard computations one can see that the Gateaux derivative of $ I $ is given by 
\begin{small}
$$
I'(u)v = \intR [|\nabla u|^{N-2}\nabla u \nabla v + |u|^{N-2}uv] dx + \intR \intR \ln(1+|x-y|) |u(x)|^{N} |u(y)|^{N-2}u(y)v(y) dx dy - \intR f(u)v dx 
$$
\end{small}
and that $ I\in C^1(X, \R) $. Also, $ I $ is clearly invariant under $ \mathbb{Z}^N $ translations. Hence, noting that any critical point of $ I $ is a weak solution for (\ref{P}), we will be concerned in finding critical points for $ I $. 

Provided the necessary background, we present our main result concerning the existence of solutions for (\ref{P}).

\bt\label{t1}
Assume $ (f_1)-(f_4) $, $ q>2N $ and that $ C_q>0 $ is sufficiently large. Then,
\begin{itemize}
\item[(i)] Equation (\ref{P}) has a solution $ u\in X\setminus \{0\} $ with $ I(u)=c_{mp} $, where 
\begin{equation}\label{eq9}
c_{mp}= \inf\limits_{\gamma \in \Gamma} \max\limits_{t\in [0, 1]} I(\gamma(t)) ,
\end{equation}
with $ \Gamma = \{ \gamma\in C([0, 1], X) \ ; \ \gamma(0)=0 \ , \ I(\gamma(t))< 0 \} $.

\item[(ii)] Equation (\ref{P}) has a non-trivial ground state solution, i.e there exists $ u\in X\setminus \{0\} $ such that $$ I(u)=c_g = \inf\{I(v) \ ; \ v\in X\setminus\{0\} \mbox{ \ is a solution of (\ref{P}}) \} . $$ 
\end{itemize}
\et

\bc\label{c400}
We have that $c_g > 0$. Particularly, problem \eqref{P} has only positive energy level non-trivial solutions.
\ec

For the second main result, we are concerned with multiplicity of solutions. However, to obtain guarantee that $ I $ has the necessary geometry, we need to ask for $ (f_5) $ and replace the condition $ (f_1) $ by the condition below.
$$ f\in C(\R , \R), \ f \mbox{ \ is odd and has subcritical exponential growth.} \leqno{(f_1 ')}$$
 
$$ \mbox{the function \ } t\mapsto \dfrac{f(t)}{t^{2N-1}} \mbox{ \ is increasing in \ } (0, + \infty) . \leqno{(f_5)}$$
From $ (f_5) $, it follows that $ \frac{f(t)}{t^{2N-1}} $ is decreasing in $ (-\infty , 0) $.

Moreover, since we can control the exponent using $ \alpha >0 $, we can change condition $ (f_4) $ by a more general condition, that is
$$\mbox{ there  exists } q>2N \mbox{ \ and \ } M_1 > 0 \ \mbox{ such that}\  F(t) \geq M_1 |t|^q \ , \ \forall \ t\in \R , \leqno{(f_4 ')}$$
 
\bt\label{t2}
Suppose $ (f_1 '), (f_2), (f_3), (f_4 '), (f_5) $. Then, problem (\ref{P}) has infinitely many solutions. 
\et

Once the functional $ I $ has the same geometry as the functional associated to the problem in \cite{[boer2]}, we will not repeat the details of the proofs. The reader can found it full detailed in \cite[Section 5]{[boer2]}. One should observe that the arguments are mainly based in the geometry of $ I $ independent of the specific form of the Laplacian type operator involved.  

Throughout this paper, we will use the following notations: $ \Ls $ denotes the usual Lebesgue space with norm $ ||\cdot ||_s $ \ ; \ $ X' $ denotes the dual space of $ X $ \ ; \ $ B_r(x) $ is the ball centred in $ x $ with radius $ r>0 $ \ ;  \ $ C, C_2, ... $ will denote different positive constants whose exact values are not essential to the exposition of arguments. 

The paper is organized as follows: in section 2 we present some technical and essential results, some of them already derived in previous works and whose application to our problem is immediate. Finally, section 3 consists in the proof of a key proposition and our main result.

\section{Framework and Preliminary Results}

In this section, we will provide some extra framework informations and present very useful technical lemmas. The results that already possess fully proofs in the other works, we will only make a fast recall, in order to focus on our new results. 

Inspired by \cite{[6]}, we start defining three auxiliary functionals, $V_1: \WumN \RA [0, \infty],$ $V_2: L^{\frac{2N^2}{2N-1}}(\RN) \RA [0, \infty) $ and $V_0: \WumN \RA \R \cup \{\infty\},$ given by 
\beqa
u \mapsto V_1(u)=\intR \intR \ln(1+|x-y|)|u(x)|^N |u(y)|^N dx dy ,
\eeqa
\beqa
u \mapsto V_2(u)=\intR \intR \ln\left(1+\dfrac{1}{|x-y|}\right)|u(x)|^N |u(y)|^N dx dy ,
\eeqa
\beqa
u \mapsto V_0(u)=V_1(u, v)-V_2(u, v)=\intR \intR \ln(|x-y|)|u(x)|^N |u(y)|^N dx dy. 
\eeqa
These definitions are understood to being over measurable functions $u: \RN \RA \R $, such that the integrals are defined in the Lebesgue sense. We also recall that, for $ r > 0 $, $ \ln(1+r) \leq r $ and
\begin{equation}\label{eq12}
\ln(1+|x-y|) \leq \ln(1+|x|)+\ln(1+|y|)  \ , \ \forall \ x, y \in \RN .
\end{equation} 

The proofs of the next two technical lemmas follow some standard arguments (see for example \cite{[6], [boer2]}).

\bp \label{p3}
\begin{itemize}
\item[(i)] The space $ X $ is compactly embedded in $ \Ls $, for all $ s\in [N, +\infty) $. 

\item[(ii)] The functionals $ V_0, V_1, V_2 $ and $ I $ are of class $ C^1 $ on $ X $. Moreover, for any $ u, v\in X $,
$$
V_{1}'(u)(v)=2N\intR \intR \ln(1+|x-y|)|u(x)|^N |u(y)|^{N-2}u(y)v(y) dxdy ,
$$
and for any $ u, v \in \LNquatro $,
$$
V_{2}'(u)(v)=2N\intR \intR \ln\left(1+\dfrac{1}{|x-y|}\right)|u(x)|^N |u(y)|^{N-2}u(y)v(y) dxdy .
$$

\item[(iii)] Functional $ V_2 $ is continuous (in fact continuously differentiable) on $\LNquatro$ .

\item[(iv)] Let $ (\un)\subset X $ such that $ \un \rightharpoonup u $ in $ X $. Then, 
$$
\lim\limits_{n\RA + \infty} \intR \intR \ln(1+|x-y|)|\un(x)|^N |u(y)|^{N-2} u(y) (\un(y) - u(y)) dx dy = 0 .
$$
\end{itemize}
\ep

From a corollary of Ergorov's Theorem, we can state the next proposition, which proof can be done similarly as \cite[Proposition 3.1]{[boer2]}, with $ p=N $.

\bp\label{p2}
Let $ u\in \LN \setminus \{0\} $, $ (\un)\subset \LN $ such that $ \un(x)\RA u(x) $ a.e in $ \RN $ and $ (\vn)\subset \LN $ bounded. Set 
$$
\omega_n = \intR \intR \ln(1+|x-y|)|\un(x)|^N |\vn(y)|^N dx dy .
$$
If $\sup\limits_{n} \omega_n < + \infty $, then $ ||\vn||_\ast $ is bounded. Moreover, if $ \omega_n \RA 0 $ and $ ||\vn||_N \RA 0 $, then $ ||\vn||_\ast \RA 0 $.
\ep
\begin{proof}
From the Ergorov's corollary, there exists $ R\in \N $, $ \delta > 0 $, $ n_0 \in \N $ and $ A\subset B_R $, such that $ A $ is measurable, $ \mu(A)>0 $ and $ \un(x) > \delta $, for all $ x\in A $ and for all $ n \geq n_0 $. Then, by some calculation, one can obtain
\beq
\omega_n \geq \dfrac{\delta^N \mu(A)}{2}  (||v||_{\ast}^{N}-\ln(1+2R)||\vn||_{N}^{N}),
\eeq
from which the result follows.
\end{proof}

Let $ u\in X $. By Hardy-Littlewood-Sobolev Inequality (HLS) \cite{[15]} with $ \alpha, \beta =0  $, $ \lambda = 1 $, $ g(x)=|u(x)|^N $, $ f(y)=|u(y)|^N $ and $ \frac{1}{q}+\frac{1}{t}+\frac{1}{N}=2 $, we have
\beq \label{v2}
|V_2(u)|\leq K_0||u||_{\frac{2N^2}{2N-1}}^{2N} \ , \ \ \forall \ u\in L^{\frac{2N^2}{2N-1}}(\RN) ,
\eeq
where $ K_0 $ is the HLS constant. Particularly, $ V_2 $ takes finite values over $ L^{\frac{2N^2}{2N-1}}(\RN)\subset \WumN $. 

\bo\label{obs8} From \eqref{eq12}, 
$$\intR \intR \ln(1+|x-y|)|u(x)|^N |v(y)|^N dxdy \leq ||u||_{\ast}^{N}||v||_{N}^{N}+ ||v||_{\ast}^{N}||u||_{N}^{N}.$$
As a consequence, $ V_1(u)\leq 2||u||_{\ast}^{N}||u||_{N}^{N} $, for all $u\in \WumN$.
\eo

Before discuss further the exponential behaviour, we would like to highlight some important inequalities obtained from our main assumptions on $ f $.

\bo\label{obs1}
Given $ \varepsilon >0 $, $ q>N $ and $ \alpha > \alpha_0 $, from $ (f_2) $ and $ (f_1) $, for all $ u\in \WumN $, there exists a constant $ b_1 >0 $ such that
\beq\label{eq1}
|F(u)|\leq \dfrac{\varepsilon}{N}|u|^N + b_1 |u|^q R(\alpha ,u) .
\eeq
Similarly, there exists constants $ b_2, b_3 >0 $ satisfying
\begin{equation}\label{eq14}
|f(u)| \leq \varepsilon |u|^{N-1} + b_2 |u|^{q-1}R(\alpha , u)
\end{equation}
\eo

Next, we provide some extremely important lemmas concerning the exponential growth. Lemma \ref{l1} is the well-known Moser-Trudinger Lemma for $ N $-Laplacian. Lemma \ref{l33} provide us with an essential inequality for potentials of $ R(\alpha , u) $ and then we have some immediate corollaries derived from both. 

\bl\label{l1}
(Moser-Trudinger Lemma \cite{[doO], [17]}) Let $ N\geq 2 $, $ \alpha >0 $ and $ u\in \WumN $. Then, 
$$
\intR [ \exp(\alpha |u|^{\frac{N}{N-1}}) - S_{N-2}(\alpha , u)] dx < \infty ,
$$
where $ S_{N-2}(\alpha, u) = \sum\limits_{k=0}^{N-2}\dfrac{\alpha^{k}}{k!}|u|^{\frac{N}{N-1}k} $. Moreover, if $ ||\nabla u||_{N}^{N} \leq 1 $, $ ||u||_N \leq M < \infty $ and $ \alpha < \alpha_N = N\omega_{N-1}^{\frac{1}{N-1}} $, where $ \omega_{N-1} $ is the $ (N-1) $-dimensional measure of $ (N-1) $-sphere, then there exists a constant $ C_0 = C(\alpha , N , M) $ such that
$$
\intR [ \exp(\alpha |u|^{\frac{N}{N-1}}) - S_{N-2}(\alpha , u)] dx \leq C(\alpha , N , M) = C_0 .
$$
\el

\bl\label{l33}
(\cite[Lemma 2.3]{[w6]}) Let $ \alpha > 0 $ and $ r>1 $. Then, for every $ \beta > r $, there exists a constant $ C_\beta = C(\beta) > 0 $ such that
$$
(\exp(\alpha |t|^{p'}) - S_{N -2}(\alpha ,t))^r \leq C_\beta (\exp(\beta \alpha |t|^{p'} - S_{N -2}(\beta \alpha , t)) ,
$$
with $ \frac{1}{p}+\frac{1}{p'}=1 $. 
\el

\bc\label{obs1}
Let $ \alpha>0 $. Then, $ R(\alpha , u)^l \in L^{1}(\RN) $, for all $ u\in \WumN $ and $ l \geq 1 $.
\ec

\bc\label{l11}
Let $ u\in \WumN $, $ r>N $, $ l\geq 1 $, $ \beta > 0 $ and $ ||u||\leq M $, for $ M>0 $ sufficiently small. Then, there exists a constant $ K_1 = K_1(\beta , N, M, l, s) > 0 $ such that
$$
\intR |u|^r R(\beta , u)^l dx \leq K_1 ||u||_{t_0}^{r} ,
$$
for some $ t_0 > N $. Moreover, there exists a constant $ K_2 = K_2(\beta , N, M, l, s) > 0 $ such that
$$
\intR |u|^r R(\beta , u)^l dx \leq K_2 ||u||_{X}^{r} .
$$
\ec

Analyzing Corollary \ref{l11}, one can see that it is possible to obtain the estimatives either controlling the norm $ ||\cdot|| $, as we presented, or controlling the exponent $ \beta $ instead. Both forms are equally important inside of the work. In the next two lemmas we guarantee that $ I $ possess the Mountain-Pass geometry. Their proof are usual and will be omitted (see \cite[Lemmas 3.4 and 3.5]{[boer2]}).

\bl\label{l6}
There exists $ \rho > 0 $ such that
\begin{equation}\label{eq6}
m_\beta = \inf \{ I(u) \ ; \ u\in X \ , \ ||u||=\beta \} \mbox{ , for all \ } 0 < \beta \leq \rho 
\end{equation}
and
\begin{equation}\label{eq7}
n_\beta = \inf \{ I'(u)(u) \ ; \ u\in X \ , \ ||u||=\beta \} \mbox{ , for all \ } 0 < \beta \leq \rho .
\end{equation}
\el
\begin{proof}
Let $ u\in X\setminus \{0\} $, $ \varepsilon > 0 $ small and $ q> N $, then
$$
I(u) \geq \dfrac{||u||^N}{N}[1-\varepsilon - C_1||u||^N - C_2||u||^{q-N}] \mbox{ \ \ \ and \ \ \ } I'(u)(u)\geq ||u||^N (1-\varepsilon - C_3||u||^N - C_4 ||u||^{q-N}).
$$
\end{proof}

\bl\label{l7}
Let $ u\in X\setminus \{0\} $, $ t>0 $ and $ q> 2N $. Then,
$$
\lim\limits_{t\RA 0} I(tu) = 0 \ \ , \ \ \sup\limits_{t>0} I(tu) < +\infty \ \ \mbox{and} \ \ I(tu)\RA + \infty \ , \ \mbox{as} \ t \RA + \infty .
$$
\el
\begin{proof}
Let $ u\in X\setminus \{0\} $, $ t>0 $ and $ q> 2N $. Then, we observe that
$$
I(tu) \leq \dfrac{t^N}{N}||u||^N+\dfrac{t^{2N}}{2N}V_0(u)-C_q t^{q}||u||_{q}^{q}. 
$$
\end{proof}

Since $ I $ has the mountain pass geometry, the mountain pass level defined in (\ref{eq9}) is well-defined and satisfies $ 0< m_\rho \leq c_{mp} < + \infty $. Moreover, there exists a Cerami sequence in the level $ c_{mp} $.

In this sense, consider a sequence $ (\un)\subset X $ satisfying
\begin{equation}\label{eq8}
0 < d = \sup\limits_{n\in \N} I(\un) < + \infty \mbox{ \ \ \ and \ \ \ } ||I'(\un)||_{X'}(1+||\un||_X) \RA 0 \mbox{ \ as \ } n \RA + \infty .
\end{equation}
Then, by a direct calculation, we get the following lemma.

\bl\label{l8}
Let $ (\un) \subset X $ satisfying (\ref{eq8}). Then, $ (\un) $ is bounded in $ \WumN $.
\el
\begin{proof}
From (\ref{eq8}) and $ (f_3) $, we have
$$
d+ o(1)  \geq I(\un) - \dfrac{1}{2N} I'(\un)(\un)\geq \dfrac{1}{2N}||\un||^N +  \left(\dfrac{\theta}{2N} - 1 \right) \intR F(\un) dx \geq \dfrac{1}{2N}||\un||^N .
$$
\end{proof}

From the above lemma, for a Cerami sequence at level $ c_{mp} $, we get $ 2Nc_{mp}+o(1) \geq ||\un||^N $. Then, we prove the next lemma similarly as \cite[Lemma 3.7]{[boer2]}, which give us a necessary way to control the exponential term for Cerami sequences in levels $ d\in (-\infty , c_{mp}] $.
\bl\label{l10}
Let $ (\un)\subset X $ be a Cerami sequence for $ I $ at level $ c_{mp} $ and $ q>2N $. Then, for some $ \rho_0 >0 $ sufficiently small,
\begin{equation}\label{eq11}
\limsup\limits_{n} ||\un|| < \rho_0 .
\end{equation}
\el
\begin{proof}
By similar arguments as in \cite[Lemma 26]{[boer2]}, we obtain 
$$
\limsup\limits_{n} ||\un||^N  \leq \dfrac{2(q-N)}{q}\dfrac{S_{q}^{\frac{qN}{q-N}}}{(qC_q)^{\frac{N}{q-N}}} ,
$$
from which the result follows for $ C_q > 0 $ sufficiently large.
\end{proof}

\section{Proof of Theorems \ref{t1} and \ref{t2}}

We are now ready to prove our main result. However, we still need a key proposition that provide us with nontrivial critical points of $ I $ in $ X $.

\bp\label{p1}
Let $ q > 2N $ and $ (\un)\subset X $ a Cerami sequence for $ I $ in level $ c_{mp} $. Then, passing to a subsequence, if necessary, only one between the alternatives occurs:

\noindent \textbf{(I)} $ ||\un||\RA 0 $ and $ I(\un)\RA 0 $.

\noindent \textbf{(II)} There exists points $ \yn \in \mathbb{Z}^N $ such that $ \yn \ast \un \RA u $ in $ X $ for a non-trivial critical point $ u\in X $ of $ I $.
\ep
\begin{proof}
From Lemmas \ref{l1}, \ref{l33}, \ref{l8}, \ref{l10}, Propositions \ref{p3} and \ref{p2}, Corollary \ref{l11}, HLS and Hölder Inequalities and the $ (S^{+}) $ property of $ N $-Laplacian, one can prove this propositon following the same steps as in the proof of \cite[Proposition 4.1]{[boer2]}.
\end{proof}

\begin{proof}[Proof of Theorem \ref{t1}]
\textbf{(i)} From Lemma \ref{l6} and Proposition \ref{p1} there exists a non-trivial critical point of $ I $, $ u_0 \in X $, such that $ I(u_0)=c_{mp} $. 

\textbf{(ii)} Define $ K=\{v \in X\setminus \{0\} \ ; \ I'(v)=0\} $. Since $ u_0\in K $, $ K\neq \emptyset $. Thus, we can consider $ (\un)\subset K $ satisfying $ I(\un)\RA c_g = \inf\limits_{v\in K} I'(v) $.

Observe that $ c_g \in [-\infty , c_{mp}] $. If $ c_g = c_{mp} $ it is done. If $ c_g < c_{mp} $, then from definition of $ K $ and Proposition \ref{p1} there exists $ (\yn)\subset \mathbb{Z}^N $ such that $ \until \RA u $ in $ X $, for a non-trivial critical point $ u $ of $ I $ in $ X $ and we conclude that $ u\in K $ and 
$
I(u)= c_g .
$
Particularly, we see that $ c_g > -\infty $.
\end{proof}

\begin{proof}[Proof of Corollary \ref{c400}]
Let $ (\until) \subset X $ be the sequence obtained in Theorem \ref{t1} such that $ \until \RA u $, for $ u \in X\setminus\{0\} $, 
$$
I(\until)=I(\un) \RA c_g \mbox{ \ \ \ and \ \ \ } ||I'(\un)||_{X'}(1+||\un||_X)\RA 0 \mbox{ , as \ } n \RA + \infty ,
$$
$ I(u)=c_g $ and $ I'(u)\equiv 0 $. Since $ ||\cdot|| $ is invariant under $ \mathbb{Z}^N$-translation and $ \until \rightharpoonup u $ in $ \WumN $, from the fact that any norm is weakly lower semicontinuous and by a similar argument as in Lemma \ref{l8}, we get that 
$$
2Nc_g + o(1) \geq ||\until||^{N} \ \ \ \Rightarrow \ \ \ 2Nc_g \geq \liminf  ||\until||^{N} \geq ||u||^N > 0 .
$$
\end{proof}

To prove multiplicity, we use the genus theory, whose definition and basic properties can be found in \cite[Chapter II.5]{struwe}. As we told, the auxiliary function $ \varphi_u : \R \RA \R $, given by $ \varphi_u(t) = I(tu) $, for all $ u\in X\setminus\{0\} $ and $ t \in \R $, enjoy the same geometry verified in \cite[Lemma 5.1]{[boer2]}, that is, 

\bl
(\cite[Lemma 5.1]{[boer2]}) \textbf{(a)} Let $ u\in X\setminus \{0\} $. Then, $ \varphi_u $ is even and there exists a unique $ t_u \in (0, +\infty) $ such that $ \varphi_u '(t)>0 $, for all $ t\in (0, t_u) $, and $ \varphi_u '(t)<0 $, for all $ t\in (t_u , \infty) $. Moreover, $ \varphi_u (t) \RA -\infty $, as $ t \RA +\infty $.

\noindent \textbf{(b)} Let $ u\in X\setminus \{0\} $. Then, there exists a unique $ t_{u}'\in (0, + \infty )$ such that $ \varphi_u(t)>0 $, for $ t\in (0, t_{u}') $, and $ \varphi_u(t)<0 $, for $ t\in (t_{u}', + \infty) $. Moreover, $ t_u $ given by item (a) is a global maximum for $ \varphi_u $.

\noindent \textbf{(c)} For each $ u\in X\setminus \{0\} $, the map $ u \mapsto t_u ' $ is continuous.
\el

As we mentioned in the introduction, we will not explicit here all the results necessaries to prove Theorem \ref{t2}, since there are many and their proof are very similar to those one in \cite[Lemma 29 - Proposition 47]{[boer2]}. So, we only define the necessary elements that appear explicit in the proof of Theorem \ref{t2}.

Consider the Nehari's manifold for $ I $, defined by
\begin{equation}\label{VN}
\Nn = \{ u\in X \setminus \{0\} \ ; \ I'(u)(u) = 0 \},
\end{equation}
the sets
$$
I^{c}=\{u\in X \ ; \ I(u)\leq c \} \mbox{ \ , for \ } c\in \R \mbox{ \ and \ } D=I^0 ,
$$
and the values
$$
c_k = \inf\{c \geq 0 \ ; \ \gamma_D(I^c)\geq k \} \mbox{ \ , \ } \forall \ n\in \N .
$$
We can see that, in the case of $ k=1 $, we have
$ c_1 = \inf\limits_{N} I = \inf\limits_{u\in X \setminus \{0\}}\sup\limits_{t > 0}I(tu) >0 $. Moreover, we have the following two propositions that are the base for the proof of Theorem \ref{t2}.

\bp\label{p112}
Let $ k\in \N $. Then, $ c_k $ is a critical value of $ I $.
\ep

\bp\label{p113}
We have $ c_k \RA + \infty $, as $ k \RA + \infty $.
\ep

\begin{proof}[Proof of Theorem \ref{t2}] 
By Proposition \ref{p113}, passing to a subsequence if necessary, $ c_k \RA + \infty $ monotonously increasing. Thus, from Proposition \ref{p112}, there exists $ \uk \in X $ satisfying $ I(\uk)=c_k $ and $ I'(\uk)=0 $, for all $ k\in \N $. Since $ (c_k) $ is monotone, $ c_k >0 $, for all $ k\geq 2 $. Then, we conclude that the functions $ u_k $ are distinct and that $ \uk \neq 0 $, for all $ n\in \N $. Moreover, as $ I $ is odd, the same holds for $ -\uk $ and $ I(\pm \uk) \RA + \infty $.  
\end{proof}

\section{Final Considerations}

In this final section, we only want to call attention the reader that all the above arguments, included those in \cite{[boer2]} can be adapted, with minor changes, for a non-constant potential $ a: \RN \RA \R $ satisfying the following conditions
$$a: \RN \RA \R \mbox{ is continuous, } \mathbb{Z}^N\mbox{-periodic , } a\in L^{\infty}(\RN) \mbox{ and } \inf\limits_{x\in \RN} a(x) = a_0 > 0 . \leqno{(a_0)}$$
One also could investigate the case where the potential $ a $ is not invariant under $ \mathbb{Z}^N $ translations but is asymptotically $ \mathbb{Z}^N $-periodic, that is, there exists a $ \mathbb{Z}^N $-periodic potential $ a_p : \RN \RA \R $ such that $ a_p $ satisfies $ (a_0) $,
$$
0 < \inf\limits_{x\in \RN} a(x) \leq a(x) \leq a_p(x) \ , \ \forall \ x \in \RN \leqno{(a_1)}
$$
and
$$
\lim\limits_{|x|\RA + \infty}|a(x) - a_p(x)| = 0 . \leqno{(a_2)}.
$$
To do that, one can argue as in \cite{[alves2]}. \\

\noindent \textbf{Acknowledgements:} The first author was supported   by  Coordination of Superior Level Staff Improvement-(CAPES)-Finance Code 001 and  S\~ao Paulo Research Foundation-(FAPESP), grant $\sharp $ 2019/22531-4, while the second  author was supported by  National Council for Scientific and Technological Development -(CNPq),   grant $\sharp $ 307061/2018-3 and FAPESP  grant $\sharp $ 2019/24901-3.

\end{document}